\documentclass[12pt]{article}
\usepackage[centertags]{amsmath}
\usepackage{amsfonts}
\usepackage{amssymb}
\usepackage{latexsym}
\usepackage{amsthm}
\usepackage{newlfont}
\usepackage{graphicx}
\usepackage{booktabs}
\date{}

\newlength{\defbaselineskip}
\newcommand{\setlinespacing}[1]%
           {\setlength{\baselineskip}{#1 \defbaselineskip}}

\newcommand{\N}{{\mathbb{N}}}

\newcommand{\actaqed}{\hfill $\actabox$}
{\medskip\noindent \textit{Proof of #1. }}%
{\actaqed \medskip}

\def \cF{\mathcal F}

\def\R{{\mathbb R}}
\def\Z{\mathbb Z}

\def \<{\langle}
\def\>{\rangle}

\def\la{\lambda}

\def\bx{\mathbf x}
\def\by{\mathbf y}
\def\bz{\mathbf z}
\def\bk{\mathbf k}
\def\bu{\mathbf u}

\def\bm{\mathbf m}

\def\bs{\mathbf s}

\newtheorem{Theorem}{Theorem}[section]
\newtheorem{Lemma}{Lemma}[section]
\newtheorem{Definition}{Definition}[section]

\numberwithin{equation}{section}

\newcommand{\be}{\begin{equation}}
\newcommand{\ee}{\end{equation}}

\begin{document}

\title{On the fixed volume discrepancy of the Fibonacci sets in the integral norms}
\author{V.N. Temlyakov\thanks{University of South Carolina, USA; Steklov Institute of Mathematics and Lomonosov Moscow State University, Russia.} \, 
and M. Ullrich\thanks{Johannes Kepler University Linz, Austria.}}
\maketitle
\begin{abstract}
This paper is devoted to the study of a discrepancy-type characteristic 
-- the fixed volume discrepancy -- 
of the Fibonacci point set in the unit square. 
It was observed recently that this new characteristic allows us to obtain 
optimal rate of dispersion from numerical integration results. 
This observation motivates us to thoroughly study this new version of 
discrepancy, which seems to be interesting by itself.  
The new ingredient of this paper 
is the use of the  average over the shifts of hat functions instead of 
taking the supremum over the shifts. We show that this change in the setting 
results in an improvement of the upper bound for the 
smooth fixed volume discrepancy, similarly to the well-known results 
for the usual $L_p$-discrepancy. 
Interestingly, this shows that ``bad boxes'' for the usual discrepancy 
cannot be ``too small''.
The known results on smooth discrepancy show that the obtained bounds 
cannot be improved in a certain sense.
\end{abstract}

\section{Introduction} 
\label{I} 

This paper is devoted to the study of a discrepancy-type characteristic -- the fixed volume discrepancy -- 
of a point set in the unit square $\Omega_2:=[0,1)^2$. 
We refer the reader to the following books and survey papers on discrepancy theory and 
numerical integration \cite{BC}, \cite{Mat}, \cite{NW10}, \cite{VTbookMA}, \cite{Bi} \cite{DTU}, \cite{T11}, and \cite{VT170}. 
 Recently, an important new observation was made in \cite{VT161}. 
It claims that a new version of discrepancy -- the $r$-smooth fixed volume discrepancy -- allows 
us to obtain optimal rate of \emph{dispersion} from numerical integration results 
(see \cite{AHR,BH19,DJ,RT,Rud,Sos,Ull,U19,UV18} for some recent results on dispersion). 
This observation motivates us to thoroughly study this new version of discrepancy, 
which seems to be interesting by itself. 

The \emph{$r$-smooth fixed volume discrepancy} takes into 
account two characteristics of a smooth hat function $h^r_B$ -- its smoothness $r$ and the volume of  its support $v:=\mathrm{vol}(B)$ (see the definition of $h^r_B$ below). The new ingredient of this paper 
is the use of the $L_p$, $1\le p<\infty$, average over the shifts of hat functions instead of 
taking the supremum over the shifts. We show that this change in the setting of the problem results in an improvement of the upper bound for the $r$-smooth fixed volume discrepancy of the special sets of points -- the Fibonacci point sets. 
For these sets with $b_n$ elements (see below), we get 
$(\log(b_n v))^{1/2}$ for $1\le p<\infty$, instead of $\log(b_n v)$ for $p=\infty$. 
The known results on $r$-smooth discrepancy show that both bounds cannot be improved in a 
certain sense (see the end of Introduction for a detailed discussion). 
The new results are only for the Fibonacci point sets, i.e., in dimension 2, 
and for $L_p$-averaging in the periodic setting, i.e., with respect to the torus geometry. 
However, we present the corresponding definitions and some known results in a general 
setting on the unit cube $\Omega_d:=[0,1)^d$. 
We now proceed to a formal description of the problem setting and to formulation of the results. 

Denote by $\chi_{[a,b)}(x)$  a univariate characteristic function (on $\R$) of the interval $[a,b)$ 
and, for $r=1,2,3,\dots$, we inductively define
$$
h^1_u(x):= \chi_{[-u/2,u/2)}(x)
$$
and 
$$
h^r_u(x) := h^{r-1}_u(x)\ast h^1_u(x),
$$
where
$$
f(x)\ast g(x) := \int_\R f(x-y)g(y)dy.
$$
Note that $h^2_u$ is the \emph{hat function}, i.e., 
$h_{u}^2(x) =\max\{u-|x|,0\}$.

Let $\Delta_tf(x):=f(x)-f(x+t)$ be the first difference. 
We say that a univariate function $f$ has smoothness $1$ in $L_1$ if 
$\|\Delta_t f\|_1\le C|t|$ for some absolute constant $C<\infty$. 
In case $\|\Delta^r_t f\|_1 \le C|t|^r$, where $\Delta^r_t:= (\Delta_t)^r$ is the $r$th difference operator, $r\in \N$, we say that $f$ has smoothness $r$ in $L_1$.
Then, $h^r_u(x)$ has smoothness $r$ in $L_1$ and has support $(-ru/2,ru/2)$. 

For a box $B$ of the form
\be\label{eq:B}
B= \prod_{j=1}^d [z_j-ru_j/2,z_j+ru_j/2)
\ee
 define
\be\label{eq:hB}
h^r_B(\bx):= h^r_\bu(\bx-\bz):=\prod_{j=1}^d h^r_{u_j}(x_j-z_j).
\ee
We begin with the \emph{non-periodic} $r$-smooth fixed volume discrepancy introduced and studied in \cite{VT161}.
\begin{Definition}\label{ED.1} Let $r\in\N$, $v\in (0,1]$ and $\xi:= \{\xi^\mu\}_{\mu=1}^m\subset [0,1)^d$ be a point set. 
We define the $r$-smooth fixed volume discrepancy with equal weights as
\be\label{E.1}
D^r(\xi,v):=  \sup_{B\subset \Omega_d:vol(B)=v}\left|\int_{\Omega_d} h_B^r(\bx)d\bx-\frac{1}{m}\sum_{\mu=1}^m h_B^r(\xi^\mu)\right|.
\ee
The optimized version of the $r$-smooth fixed volume discrepancy is defined as follows
\be\label{E.2}
D^{r,o}(\xi,v):=  \inf_{\la_1,\dots,\la_m}\sup_{B\subset \Omega_d:vol(B)=v}\left|\int_{\Omega_d} h_B^r(\bx)d\bx- \sum_{\mu=1}^m \la_\mu h_B^r(\xi^\mu)\right|.
\ee
\end{Definition}
\noindent
Clearly, we have $D^{r,o}(\xi,v) \le D^{r}(\xi,v)$.

It is well known that the Fibonacci cubature formulas are optimal in the sense 
of order for numerical integration of different kind of smoothness classes of 
functions of two variables, see e.g.~\cite{DTU,TBook,VTbookMA}. 
We present a result from \cite{VT161}, which shows that the Fibonacci point set 
has good fixed volume discrepancy.

Let $\{b_n\}_{n=0}^{\infty}$, $b_0=b_1 =1$, $b_n = b_{n-1}+b_{n-2}$,
$n\ge 2$, be the Fibonacci numbers. Denote the $n$th {\it Fibonacci point set} by
$$
\cF_n:= \Big\{\big(\mu/b_n,\{\mu b_{n-1} /b_n \}\big)\colon\, \mu=1,\dots,b_n\Big\}.
$$
In this definition $\{a\}$ is the fractional part of the number $a$.   The cardinality of the set $\cF_n$ is equal to $b_n$. In \cite{VT161} we proved 
the following upper bound.

\begin{Theorem}\label{ET.1} Let $r\ge2$. 
There exist constants $c,C>0$ such that for any $v\ge c/b_n$ 
we have 
\be\label{E.3}
D^{r}(\cF_n,v)
\,\le\, C\,\frac{\log(b_n v)}{b_n^r}.
\ee
\end{Theorem}

\bigskip

The main object of our interest in this paper is the \emph{periodic} $r$-smooth $L_p$-discrepancy 
of the Fibonacci point sets. 
For this, we define the periodization $\tilde f$ 
(with period $1$ in each variable) 
of a function $f\in L_1(\R^d)$ with a compact support by
$$
\tilde f(\bx) := \sum_{\bm\in \Z^d} f(\bm+\bx)
$$
and, for each $B\subset [0,1)^d$, 
we let $\tilde h^r_B$ be the periodization of $h^r_B$ from~\eqref{eq:hB}.

We now define the periodic $r$-smooth $L_p$-discrepancy.
\begin{Definition}\label{ID.1}For $r\in\N$, $1\le p\le \infty$ and $v\in(0,1]$ define the 
periodic $r$-smooth fixed volume $L_p$-discrepancy of a point set $\xi$ by  
\be\label{1.4}
  \tilde D^{r}_p(\xi,v):=
	\sup_{B\subset \Omega_d:vol(B)=v}
	\left\|\int_{\Omega_d} \tilde h^r_B(\bx-\bz)d\bx- \frac1m\sum_{\mu=1}^m \tilde h^r_B(\xi^\mu-\bz)\right\|_p
\ee
where the $L_p$-norm is taken with respect to $\bz$ over the unit cube $\Omega_d=[0,1)^d$.\\
Analogously to \eqref{E.2} we may define the optimized version $\tilde D^{r,o}_p(\xi,v)$.
\end{Definition}

In the case of $p=\infty$ this concept was introduced and studied in \cite{VT163}. 


We prove the following upper bound for $1\le p<\infty$.
\begin{Theorem}\label{IT.1} 
Let $r\in\N$ and $1\le p<\infty$. 
There exist constants $c,C>0$ such that for any $v\ge c/b_n$ 
we have 
$$
\tilde D^{r}_p(\cF_n,v) \,\le\, C\,\frac{\sqrt{\log(b_n v)}}{b_n^r}.
$$
\end{Theorem}

\medskip

In  the case $p=\infty$ we prove a weaker upper bound.

\begin{Theorem}\label{IT.2} Let $r\ge 2$. There exist constants $c,C>0$ such that for any $v\ge c/b_n$ 
we have 
$$
\tilde D^{r}_\infty(\cF_n,v) \,\le\, C\,\frac{\log(b_n v)}{b_n^r}.
$$
\end{Theorem}


We now give some comments, which show that Theorems \ref{IT.1} and \ref{IT.2} cannot be improved in a certain sense. We do not know if Theorems \ref{IT.1} and \ref{IT.2} are sharp 
for all $v\ge c/b_n$. 
The known results show that these theorems are sharp in some cases for the 
supremum over $v$. The following quantities
$$
  \tilde D^{r,o}_{p,\infty}(\xi):= \sup_{v}\tilde D^{r,o}_p(\xi,v)
$$
have been studied in \cite{VT163} and \cite{VT164}. 
We cite some results from there. The following lower bound follows from 
stronger results in \cite{VT164}.
Let $r,d\in\N$. Then for any point set $\xi\subset\Omega_d$ with $\#\xi=m$ we have
\be\label{1.5}
    \tilde D^{r,o}_{2,\infty}(\xi) \geq   C(r,d)\,   m^{-r}(\log   m)^{(d-1)/2},   \qquad
C(r,d)>0.
\ee
Under an extra assumption on $r$, namely, assuming that $r$ is an even number, 
we can derive an extended to $p>1$ inequality~\eqref{1.5} 
for $\tilde D^{r,o}_{p,\infty}(\xi)$ from \cite{VT164}. 
In the case $r=1$ the quantity under consideration corresponds to 
the classical (non-smooth) discrepancy, and the above mentioned bounds were already 
proven in~\cite{Ro54,Sch77}.

For $p=\infty$ the following result was proved in \cite{VT163}. 
For any point set $\xi\subset\Omega_d$ with $\#\xi=m$ 
we have for even integers $r$ that
\be\label{eq:lb-discr}
\tilde D^r_\infty(\xi) \ge C(r,d)\, m^{-r}(\log m)^{d-1}
\ee
with a positive constant $C(r,d)$. 
This result even holds if we allow weights (as in the optimized version) 
$\lambda_1,\dots,\lambda_m$ satisfying condition  
$$
\sum_{\mu=1}^m|\lambda_\mu| \le B
$$
for some fixed $B<\infty$.

Finally, let us add that Theorems \ref{IT.1} and \ref{IT.2} show that 
the ``bad boxes'', i.e., the boxes that fulfill the lower bounds~\eqref{1.5} 
or~\eqref{eq:lb-discr}, 
must have volume at least $m^{-1+\delta}$ for some fixed $\delta>0$.
This is interesting as one might think that boxes of volume at most 
$(\log m)^c/m$ (for some large $c$) may already suffice.

\section{Proofs of Theorems \ref{IT.1} and \ref{IT.2}}
\label{Fib}

The proofs of both theorems go along the same lines. We give a detailed proof of Theorem \ref{IT.1} and point out a change of this proof, which gives Theorem \ref{IT.2}.
For continuous functions of two
variables, which are $1$-periodic in each variable, define
cubature formulas
$$
\Phi_n(f) :=b_n^{-1}\sum_{\mu=1}^{b_n}f\bigl(\mu/b_n,
\{\mu b_{n-1} /b_n \}\bigr),
$$
 called the {\it Fibonacci cubature formulas}.   Denote 
$$
\mathbf y^{\mu}:=\bigl(\mu/b_n, \{\mu
b_{n-1}/b_n\}\bigr), \quad \mu = 1,\dots,b_n,
$$
and
$$
\Phi(\mathbf k) \,:=\, \Phi_n\left(e^{i2\pi(\mathbf k,\mathbf x)}\right)
\,=\, b_n^{-1}\sum_{\mu=1}^{b_n}e^{i2\pi(\mathbf k,\mathbf y^{\mu})}.
$$
Note that
\be\label{2.1}
\Phi_n (f) =\sum_{\mathbf k}\hat f(\mathbf k)\, \Phi(\mathbf k),\quad 
\hat f(\bk) := \int_{[0,1)^2} f(\bx)\, e^{-i2\pi(\bk,\bx)}d\bx,
\ee
where for the sake of simplicity we may assume that $f$ is a
trigonometric polynomial. It is clear that (\ref{2.1}) holds for $f$ with absolutely convergent Fourier series. 

It is easy to see that the following relation holds
\be\label{2.2}
\Phi(\mathbf k)=
\begin{cases}
1&\quad\text{ for }\quad \mathbf k\in L(n),\\
0&\quad\text{ for }\quad \mathbf k\notin L(n),
\end{cases}
\ee
where
$$
L(n) :=\Bigl\{ \mathbf k = (k_1,k_2)\in\Z^2\colon\; k_1 + b_{n-1} k_2\equiv 0
\; \pmod {b_n}\Bigr\}.
$$
  For $N\in\N$ define the {\it hyperbolic cross} (in dimension~2) by
$$
\Gamma(N):= \left\{\bk= (k_1,k_2)\in\Z^2\colon \prod_{j=1}^2 \max(|k_j|,1) \le N\right\}.
$$
The following lemma is well known (see, for instance, \cite{VTbookMA}, p.274).

\begin{Lemma}\label{L2.1} There exists an absolute constant $\gamma > 0$
such that for any $n > 2$ 
we have
$$
\Gamma(\gamma b_n)\cap \bigl(L(n)\backslash\{\mathbf 0\}\bigr) =
\varnothing.
$$
\end{Lemma}

\bigskip

Considering our (univariate) \emph{test functions} $h^r_u$ we obtain
by the properties of convolution that
$$
\hat h^r_u(y)  = \hat h^{r-1}_u(y)\, \hat h^1_u(y),\qquad y\in\R
$$
which implies for $y\neq 0$
$$
\hat h^r_u(y) = \left(\frac{\sin(\pi yu)}{\pi y}\right)^r.
$$
Therefore,
$$
\left|\hat {\tilde h}^r_u(k)\right| 
\,\le\, \min\left(|u|^r,\frac{1}{|k|^r}\right) 
\,=\, \left(\frac{|u|}{k'}\right)^{r/2}\min\left(|k' u|^{r/2},\frac{1}{|ku|^{r/2}}\right), 
$$
where $k':=\max\{1,|k|\}$. 
(Here, we used for a moment $\hat h$ for the Fourier transform of $h$ on $\R$. 
This should not lead to any confusion.)

\bigskip

We now proceed with some considerations in arbitrary dimension~$d$. 
It is convenient for us to use the following abbreviated notation for the product
$$
pr(\bu):= pr(\bu,d) := \prod_{j=1}^d u_j.
$$
For $B\subset\Omega_d$ of the form~\eqref{eq:B} and $\bz\in\Omega_d$, we have 
\be\label{eq:shift}
\hat {\tilde h}^r_{B+\bz}(\bk) \,=\, e^{-i2\pi(\bk,\bz)}\, \hat {\tilde h}^r_B(\bk),
\ee
where $\tilde h^r_{B+\bz}(\bx):=\tilde h^r_{B}(\bx-\bz)$, see~\eqref{eq:hB}. 
Therefore, we obtain from the above that
\[
\left|\hat {\tilde h}^r_B(\bk)\right| 
\,\le\, \prod_{j=1}^d \left(\frac{|u_j|}{k_j'}\right)^{r/2}\min\left(|k_j' u_j|^{r/2},\frac{1}{|k_j u_j|^{r/2}}\right).
\]
For $\bs\in \N_0^d$, 
we define
$$
\rho(\bs) := \Big\{\bk \in \Z^d\colon [2^{s_j-1}] \le |k_j| < 2^{s_j}, 
\quad j=1,\dots,d \Big\},
$$
where $[a]$ denotes the integer part of $a$, and obtain, 
for $\bk\in\rho(\bs)$, that  
\be\label{eq:HB}
\left|\hat {\tilde h}^r_B(\bk)\right| 
\,\le\, H_B^r(\bs) 
\,:=\, \left(\frac{pr(\bu)}{2^{\|\bs\|_1}}\right)^{r/2}\,\prod_{j=1}^d \min\left((2^{s_j}u_j)^{r/2},\frac{1}{(2^{s_j}u_j)^{r/2}}\right).
\ee
Later we will need certain sums of these quantities. 
First, consider
$$
\sigma^r_\bu(t) \,:=\, \sum_{\|\bs\|_1=t}\prod_{j=1}^d \min\left((2^{s_j}u_j)^{r/2},\frac{1}{(2^{s_j}u_j)^{r/2}}\right),\quad t\in\N_0.
$$
The following technical lemma is part~(I) from \cite[Lemma 6.1]{VT161}. 

\begin{Lemma}\label{L2.2} Let $r>0$, $t\in \N$ and $\bu\in (0,1/2]^d$ 
be such that $pr(\bu)\ge 2^{-t}$. Then, we have 


\[
\sigma^r_\bu(t) \,\le\, C(d) \frac{\left(\log(2^{t+1}pr(\bu))\right)^{d-1}}{(2^tpr(\bu))^{r/2}}.
\]
%
\end{Lemma}

\bigskip

\noindent
This lemma and \eqref{eq:HB} imply that
\be\label{eq:HB-bound}
\sum_{\|\bs\|_1=t} H_B^r(\bs)^2 
\,\le\, C_1\, 2^{-2rt}\, \Big(\log\big(2^{t}\, v\big)\Big)^{d-1}, 
\ee
where $v:=\mathrm{vol}(B)=r^d\, pr(\bu)$, 
for all $r\ge1$ and
all $t\in\N_0$ with $v\ge r^d\,2^{-t+1}$
and an absolute constant $C_1<\infty$.

\bigskip

Additionally, we need a result from harmonic analysis -- a corollary of the Littlewood-Paley theorem. Denote
$$
\delta_\bs(f,\bx):= \sum_{\bk\in\rho(\bs)} \hat f(\bk)e^{i2\pi(\bk,\bx)}.
$$
Then it is known that for $p\in [2,\infty)$ one has
\be\label{eq:LP}
\|f\|_p \le C(d,p) \left(\sum_{\bs\in\N_0^d}\|\delta_\bs(f)\|_p^2\right)^{1/2}.
\ee
Note that in the proof of Theorem \ref{IT.2} we use 
the simple triangle inequality
\be\label{eq:triangle}
\|f\|_\infty \le  \sum_{\bs}\|\delta_\bs(f)\|_\infty.
\ee
instead of \eqref{eq:LP}.

\bigskip

We are now considering the case $d=2$.
Let us define
\be\label{2.3}
E^r_B(\bz):=  \frac{1}{b_n}\sum_{\mu=1}^{b_n}  \tilde h^r_B(\by^\mu-\bz)- \int_{[0,1)^d} \tilde h^r_B(\bx)d\bx 
\ee 
such that 
\[
\tilde D^{r}_p(\cF_n,v) \,=\, \sup_{B\subset \Omega_d:vol(B)=v}
	\left\|E^r_B\right\|_p
\]
By formulas \eqref{2.1}, \eqref{2.2} and \eqref{eq:shift} we obtain
\[
E^r_B(\bz) \,=\, \sum_{\bk\neq 0} \hat {\tilde h}^r_B(\bk)\,\Phi(\bk)\, e^{-i2\pi(\bk,\bz)}
\,=\, \sum_{\bk\in L(n)\setminus\{0\}} \hat {\tilde h}^r_B(\bk)\, e^{-i2\pi(\bk,\bz)}.
\]
It is apparent from~\eqref{eq:LP} that it remains to bound 
$\|\delta_s(E_B^r)\|_p$. 

If $t\neq 0$ is such that $2^t\le \gamma b_n$ then for
$\bs$ with $\|\bs\|_1=t$ we have $\rho(\bs)\subset \Gamma(\gamma b_n)$. 
Lemma \ref{L2.1} then implies that 
$\Phi(\bk)=0$ for $\bk\in\rho(\bs)$ and, therefore, 
$\delta_s(E^r_B)=0$. 
Let $t_0\in \N$ be the smallest number satisfying 
$2^{t_0}>\gamma b_n$, 
i.e., $t_0\ge\log(b_n)-c$ for some $c<\infty$. 
Then, from~\eqref{eq:LP} for $p\in[2,\infty)$, we have 
\be\label{eq:EB-bound}
\|E^r_B\|_p \le C(p) \left(\sum_{t=t_0}^\infty \sum_{\|\bs\|_1=t}\|\delta_s(E^r_B)\|_p^2\right)^{1/2}.
\ee
Moreover, Lemma \ref{L2.1} implies that for $t\ge t_0$ we have
\be\label{eq:number}
\#\big(\rho(\bs)\cap L(n)\big) \le C_2\, 2^{t-t_0}, \quad \|\bs\|_1=t.
\ee
By Parselval's identity we obtain
\[
\|\delta_s(E^r_B)\|_2 \,=\, \sqrt{\sum_{\bk\in\rho(\bs)\cap L(n)}|\hat {\tilde h}^r_B(\bk)|^2}
\,\le\, \sqrt{\#\big(\rho(\bs)\cap L(n)\big)}\cdot H_B^r(\bs)
\]
and, by the triangle inequality, 
\[
\|\delta_s(E_B)\|_\infty \,\le\, \#\big(\rho(\bs)\cap L(n)\big)\cdot H_B^r(\bs)
\]
Hence, using the inequality 
$$
\|f\|_p \le \|f\|_2^{2/p}\|f\|_\infty^{1-2/p}
$$
for $2\le p\le \infty$, we get 
\[
\|\delta_s(E^r_B)\|_p \,\le\, \Big(\#\big(\rho(\bs)\cap L(n)\big)\Big)^{1-1/p}\cdot H_B^r(\bs).
\]
Combining this with \eqref{eq:HB-bound} for $d=2$, \eqref{eq:EB-bound} and \eqref{eq:number}, 
we finally obtain for all $v=\mathrm{vol}(B)\ge 2r^d 2^{-t_0}$ and $p\in[2,\infty)$ that
\[\begin{split}
\|E^r_B\|_p \,&\le\, C \left(\sum_{t=t_0}^\infty 2^{2(t-t_0)(1-1/p)}  
	\sum_{\|\bs\|_1=t}H_B^r(\bs)^2\right)^{1/2} \\
\,&\le\, C' \left(\sum_{t=t_0}^\infty 2^{2(t-t_0)(1-1/p)}\, 
	2^{-2rt}\, \log\big(2^{t}\, v\big)\right)^{1/2} \\
\,&=\, C'\, 2^{-r t_0} \left(\sum_{t=0}^\infty 2^{2t(1-1/p-r)}\, 
	\log\big(2^{t+t_0}\, v\big)\right)^{1/2} \\
\,&\le\, C''\, 2^{-r t_0}\, \sqrt{\log\big(2^{t_0}\, v\big)} 
	\left(\sum_{t=0}^\infty t\, 2^{2t(1-1/p-r)}\right)^{1/2}.
\end{split}\]
Using $t_0\ge\log(b_n)-c$ and that $\|E^r_B\|_p\le\|E^r_B\|_2$ for $p<2$, 
this implies Theorem~\ref{IT.1}. 
(Here, we used that clearly $r> 1-1/p$ for $p<\infty$.)

As we pointed out above, in the proof of Theorem \ref{IT.2} 
we use inequality \eqref{eq:triangle} instead of \eqref{eq:LP}. 
Moreover we use 
\[
\sum_{\|\bs\|_1=t} H_B^r(\bs) 
\,\le\, C_1\, 2^{-rt}\,\log\big(2^{t}\, v\big), 
\]
for all $r\ge1$ instead of~\eqref{eq:HB-bound}. 
However, note that we need $r>1$ for the last series in the 
above computation to be finite.
This implies
$$
\|E_B^r\|_\infty \,\le\, C\, b_n^{-r}\log\left(b_n\,v\right). 
$$

{\bf Acknowledgment.}   The work was supported by the Russian Federation Government Grant No. 14.W03.31.0031.

\end{document}